\documentclass[10pt,a4paper,reqno]{amsart}

\usepackage{latexsym}
\usepackage{verbatim}
\usepackage{epsfig}
\usepackage{rotating}
\usepackage{amssymb}
\usepackage[latin1]{inputenc}
\usepackage{afterpage}
\usepackage{color}
\usepackage{dsfont}

\usepackage{amssymb,amsfonts,amsmath,stmaryrd}

\graphicspath{{Figures/}}

\def\cE{\mathcal{E}}
\def\cR{\mathcal{R}}

\catcode`\@=11
\def\section{\@startsection{section}{1}%
 \z@{.7\linespacing\@plus\linespacing}{.5\linespacing}%
 {\normalfont\bfseries\scshape\centering}}

\def\subsection{\@startsection{subsection}{2}%
  \z@{.5\linespacing\@plus\linespacing}{.5\linespacing}%
  {\normalfont\bfseries\scshape}}

\def\subsubsection{\@startsection{subsubsection}{3}%
 \z@{.5\linespacing\@plus\linespacing}{-.5em}
  {\normalfont\bfseries\itshape}}
\catcode`\@=12

\newtheorem{Theorem}{Theorem}
\newtheorem{Lemma}[Theorem]{Lemma}

\newtheorem{claim}[Theorem]{Claim}

\newfont{\bbold}{msbm10 scaled \magstep1}
\newfont{\bbolds}{msbm7 scaled \magstep1}

\newcommand{\cT}{\mathcal T}

\newcommand{\cS}{\mathcal S}

\newcommand{\cU}{\mathcal U}
\newcommand{\cV}{\mathcal V}

\newcommand{\beq}{\begin{equation}}
\newcommand{\eeq}{\end{equation}}

\def\emm#1,{{\em #1}}

\newcommand{\cD}{\mathcal D}
\newcommand{\cF}{\mathcal F}

\def\cM{\mathcal{M}}
\def\cI{\mathcal{I}}

 \def\cB{\mathcal{B}}
 \def\cO{\mathcal{O}}

\begin{document}
\title
[Bijective counting of involutive Baxter permutations]
{Bijective counting of involutive Baxter permutations}

\author[\'E. Fusy]{\'Eric Fusy}
\thanks{Supported by the European Research Council under
the European Community's 7th Framework Programme, ERC grant agreement no 208471 - ExploreMaps project}
\address{\'E. Fusy: LIX, \'Ecole Polytechnique, 91128 Palaiseau Cedex, France}
\email{fusy@lix.polytechnique.fr}



\begin{abstract}
We enumerate bijectively the family of involutive Baxter permutations
according to various parameters; in particular we obtain an elementary proof
that the number of involutive Baxter permutations of size $2n$ with
no fixed points is $\frac{3\cdot 2^{n-1}}{(n+1)(n+2)}\binom{2n}{n}$,
a formula originally discovered by M. Bousquet-M\'elou using generating functions.
The same coefficient also enumerates planar maps with $n$ edges, endowed with an acyclic orientation having a unique
source, and such that the source and sinks are all incident to the outer face.
\end{abstract}
\maketitle

\date{\today}


\section{Introduction}
Baxter permutations, named after Glen Baxter~\cite{Bax:64} who introduced them in an analysis context, are pattern-avoiding permutations 
(precisely the forbidden patterns are $2-41-3$ and $3-14-2$) 
with many nice combinatorial properties~\cite{viennot-baxter,CGHK:78,Mal:79,dulucq-guibert-baxter}. 
Their counting coefficients appear recurrently in combinatorics;
 the so-called \emph{Baxter number} (the number of Baxter permutations of size $n$) 
$$B_n=\frac{2}{n(n+1)^2}\sum_{r=0}^{n-1}\binom{n+1}{r}\binom{n+1}{r+1}\binom{n+1}{r+2}$$ 
also counts plane bipolar orientations with $n$ edges~\cite{Bax:01},
certain rectangulations with $n$ points on the diagonal~\cite{felsner-baxter, ackerman-floorplans}, 
certain Young tableaux with
parity constraints~\cite{TS:10}, and so on. Subfamilies of Baxter permutations 
have also been considered: 
 \emph{alternating} and \emph{doubly alternating} 
  Baxter permutations have been enumerated in~\cite{CDV:86,DG:96,GL:00}, 
  and Baxter 
 permutations of size $n$ avoiding the pattern $2-4-1-3$ have been shown to be 
 in bijection with rooted non-separable maps with $n+1$ edges~\cite{DGW:96,BoBoFu09} (therefore there are 
$\tfrac{2(3n+3)!}{(2n+3)!(n+2)!}$ such permutations of size
 $n$). 
A permutation in $\frak{S}_n$ is classically drawn as an $n\times n$ grid $G$ of unit squares, with exactly
one boxed square in each row and in each column. 
With this representation in mind, 
it is known (see e.g.~\cite{BoBoFu09}) that the set of Baxter permutations of size $n$
is \emph{globally invariant} by any of the $8$ transformations of the dihedral group acting on the grid $G$.
So it is a natural problem to try to count how many Baxter permutations are fixed by a given transformation of the dihedral group. 
In a previous paper~\cite{felsner-baxter} the case of the half-turn rotation was solved (whereas the case of rotations of order $4$ is open);
the idea is that a Baxter permutation can be encoded by a nonintersecting triple of paths, in a way that commutes with the half-turn transformation.

In this note we focus on the mirror reflection according to a diagonal, that is, we count \emph{involutive} Baxter permutations.  
A difficulty is that the encoding of Baxter permutations by triples of paths does not commute in any sense with the diagonal mirror transformation
(there is no nice transfer of mirror symmetry from the Baxter permutation to the associated triple of paths).    
 An important ingredient here 
 is the recent article~\cite{BoBoFu09}, in which the authors 
 establish a direct bijective correspondence 
between Baxter permutations and so-called \emph{plane bipolar orientations},  
which are acyclic orientations on embedded planar graphs (i.e., planar maps) 
 with a unique source and
a unique sink both lying in the outer face. 
Thanks to this correspondence, the successive combinatorial manipulations
to encode involutive Baxter permutations 
 can be
carried out on (oriented) planar maps, which we find convenient
to handle due to their more geometric flavor.
As a consequence of our bijective encoding  we obtain:
\begin{itemize}
\item
 a closed-form  
 multivariate formula (in Theorem~\ref{theo:multi}) for the number of involutive Baxter permutations 
 according to the numbers of elements, descents (which are of two types, either crossing or not crossing the diagonal
$\{x=y\}$ in the diagrammatic representation), and fixed points. 
\item
a closed-form univariate formula (in Theorem~\ref{theo:uni}) 
 for the 
number $b_n$ of  involutive fixed-point free Baxter permutations
of $2n$ elements: 
$$b_n=\frac{3\cdot 2^{n-1}}{(n+1)(n+2)}\binom{2n}{n}.$$
\end{itemize}
Note that $b_n$ has surprisingly a simpler (summation-free) expression than the Baxter number $B_n$. 
The univariate formula for $b_n$ (already announced in~\cite{BoBoFu09}) 
and a multivariate formula restricted to fixed-point free Baxter permutations 
have been discovered  by M. Bousquet-M\'elou~\cite{bo10} 
using generating functions and the so-called ``obstinate'' kernel method.
As follows from the correspondence with plane orientations (to be described in Section~\ref{sec:baxter_decorated_plane}) 
the number $b_n$ also counts acyclic orientations on planar maps with $n$ edges, a unique source, and all 
extremal vertices (source and sinks) lying in the outer face. 

\vspace{.2cm}

\noindent\emph{Outline.} The main steps of our method are the following:
(i) by a quotient-argument already outlined in~\cite{BoBoFu09},   
interpret the plane bipolar orientations corresponding to involutive Baxter 
permutations as certain plane bipolar orientations with decorations
at the corners and edges incident to the sink, (ii) adapt the known
bijective encoding of plane bipolar orientations by non-intersecting triples of paths
to take account of the decorations, (iii) count the obtained non-intersecting triples
of paths using the Lindstr\"om-Gessel-Viennot lemma (for the multivariate formula)
or similar principles with small adjustments (for the univariate formula).

\section{Involutive Baxter permutations as decorated  plane  bipolar orientations}\label{sec:baxter_decorated_plane}
\begin{figure}[htb!]
\begin{center}
\includegraphics[width=13cm]{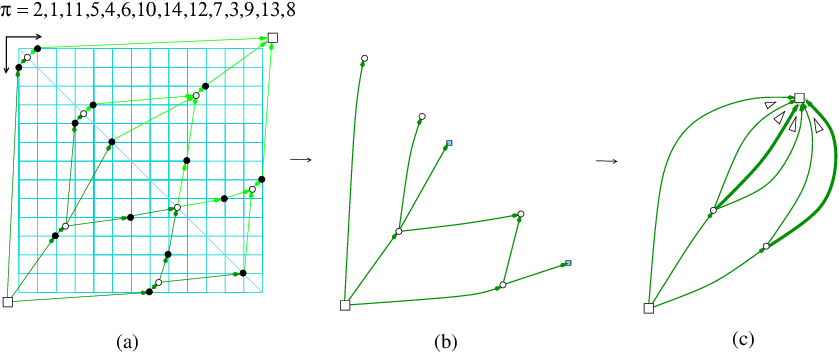}
\caption{(a) An involutive Baxter permutation $\pi$ (black points give the diagrammatic representation
of $\pi$) superimposed with the associated embedded plane bipolar orientation $O$ (a black
point in each edge of $O$, a white vertex for each vertex of $O$). The whole drawing
is invariant by the reflexion according to the line $\{x=y\}$. (b) The monosource
orientation $M$ obtained by keeping the part of $O$ below the line $\{x=y\}$, where sailing
sinks are represented as small shaded squares, (c) the associated decorated plane bipolar
orientation obtained by merging all the sinks of $M$. The marked sink-edges (which are drawn bolder)   
are those incident to sailing sinks in $M$, the marked sink-corners (indicated
by the small triangles) are the newly created
sink-corners.}
\label{fig:inv}
\end{center}
\end{figure}

Let $\cB$ be the class of Baxter permutations and $\cO$ the class of plane bipolar orientations.
As shown in~\cite{BoBoFu09} (see Figure~\ref{fig:inv} for an example), one can construct from the diagrammatic
representation of $\pi\in\cB$ an embedded plane bipolar orientation $\phi(\pi)$ where black points 
of degree 2 correspond to edges (one such vertex on each edge) and white points correspond
to vertices. The induced mapping $\Phi$ from $\cB$ to $\cO$ (where $\Phi(\pi)$ is the plane bipolar
orientation induced by $\phi(\pi)$ after erasing the black vertices) is a bijection that satisfies
several parameter correspondences (elements are mapped to edges, 
descents are mapped to non-pole vertices,...) and preserves many symmetries;
in particular if $\pi$ is involutive (i.e., $\pi^{-1}=\pi$) then $\phi(\pi)$ is fixed
by the reflection according to the line $\{x=y\}$, see Figure~\ref{fig:inv}(a). 

A \emph{planar map} is a connected graph planarly embedded in the plane (considered up to continuous deformation). 
Define a \emph{monosource} orientation as an 
acyclic orientation $O$ of a planar map with a 
unique source and such that the source and all sinks lie in the outer face. 
Additionally an arbitrary subset of the sinks
of degree $1$ are marked; these are called the \emph{sailing sinks} of $O$ 
(small shaded squares in Figure~\ref{fig:inv}(b)). 
Edges incident to sailing sinks are called \emph{sailing edges}. 

As illustrated in Figure~\ref{fig:inv}(a)-(b),  
there is a bijection between the class $\cI$ of involutive Baxter permutations and 
the class $\cM$ of monosource orientations 
which transforms standard parameters as follows (a \emph{descent} 
in a permutation $\pi\in\frak{S}_n$ is an integer $i\in[1..n-1]$ such that $\pi(i)>\pi(i+1)$,
and a descent is said to \emph{cross the diagonal} if $\pi(i)>i$ and $\pi(i+1)<i+1$):
\begin{center}
\begin{tabular}{rcl}
 $2n$ non-fixed points  & $\leftrightarrow$ &  $n$ non-sailing edges,\\ 
 $2k$ descents not crossing the diagonal & $\leftrightarrow$ &  $k$ non-extremal vertices,\\
 $p$ fixed points & $\leftrightarrow$ & $p$ sailing sinks,\\
 $r$ descents crossing the diagonal & $\leftrightarrow$ & $r$ non-sailing sinks.
 \end{tabular}
 \end{center}

We now claim that orientations in $\cM$ correspond to plane bipolar orientations
with certain decorations. Given a plane bipolar orientation, the \emph{sink-degree} is the degree of the sink, 
and a \emph{sink-edge} is an edge incident to the sink.
A \emph{corner} of a planar map is an angular sector delimited by two consecutive edges around a vertex. 
For a plane bipolar orientation, a \emph{sink-corner} is a corner incident to the sink but not in the outer face 
(note that the number of sink-corners is the sink-degree minus $1$). 
Define a \emph{decorated plane bipolar orientation} as a plane bipolar orientation
where an arbitrary subset of the 
sink-corners are marked, and a subset of sink-edges are marked in such a way
that the sink-corners incident to a marked sink-edge are marked. 
As shown in Figure~\ref{fig:inv}(b)-(c),   
there is a bijection between the class $\cM$ of monosource orientations
and the class $\cD$ of decorated plane bipolar orientations with the following
parameter-correspondence:
\begin{center}
\begin{tabular}{rcl}
 $n$ non-sailing edges  & $\leftrightarrow$ &  $n$ non-marked edges,\\ 
 $k$ non-extremal vertices & $\leftrightarrow$ &  $k$ non-pole vertices,\\
 $p$ sailing sinks & $\leftrightarrow$ & $p$ marked sink-edges,\\
 $r$ non-sailing sinks & $\leftrightarrow$ & $p+r-1$ marked sink-corners.
 \end{tabular}
 \end{center}

\section{Encoding by paths}

\begin{figure}[htb!]
\begin{center}
\includegraphics[width=\linewidth]{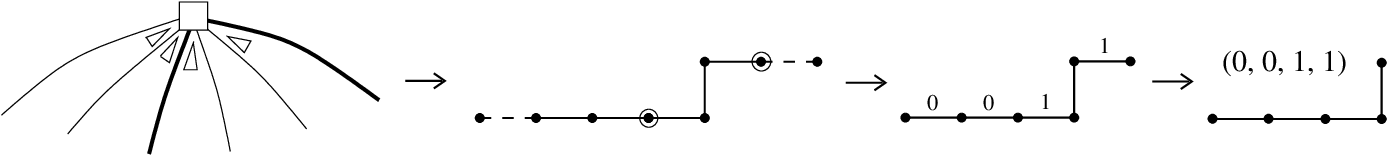}
\caption{Encoding the marked sink-corners (indicated as small triangles) 
and sink-edges (indicated as bolder edges) by a binary walk and a sequence of weights.}
\label{fig:sink}
\end{center}
\end{figure}

We now explain how to encode the sink-edges and sink-corners of a decorated 
plane bipolar orientation $O$, this is illustrated in Figure~\ref{fig:sink}.
Let $i+1$ be the sink-degree (so there are $i$ sink-corners), $p$ the number of marked
sink-edges and $q$ the number of marked sink-corners. 
By a \emph{binary walk} we will mean an oriented walk in $\mathbb{Z}^2$ having steps
 East $(+1,0)$ and North $(0,+1)$. First, encode the marked sink-corners by a binary walk 
obtained by reading the sink-corners from left to right, writing an East-step
if the sink-corner is marked and a North-step otherwise. Then append an East-step
to both the beginning and end of the binary walk; we thus have a binary walk 
with $q+2$ East-steps and $i-q$ North-steps. Note that the sink-edges that are
allowed to be marked correspond 
to points of the walk preceded and followed by an East-step (the ones corresponding
to marked sink-edges are surrounded in Figure~\ref{fig:sink}). Delete these vertices, 
and then renormalize the path to have only steps of length $1$. Above each East-step of the 
renormalized path write the number of points that have been deleted in the 
corresponding horizontal portion (before renormalization); 
finally delete the last (East-) step of the walk. 
The finally obtained walk $W$ starts with an East-step if not empty (it is empty iff all sink-edges are marked);
it has length $i+1-p$, $i-q$ North-steps, $q+1-p$ East-steps, 
and is accompanied by a sequence $S$ of $q+2-p$ nonnegative numbers adding up to $p$. 
The pair $(W,S)$ is called the \emph{sink-code} for $O$.

\begin{figure}
\begin{center}
\includegraphics[width=11.4cm]{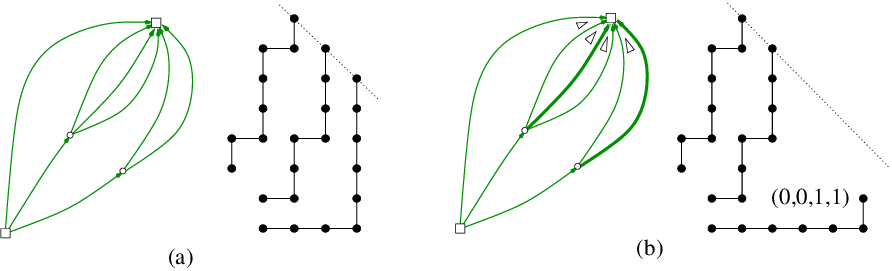}
\end{center}
\caption{(a) A plane bipolar orientation is encoded by a non-intersecting triple of binary walks. (b) The encoding
of a decorated plane bipolar orientation.}
\label{fig:bije}
\end{figure}

Let $\cT$ be the class of non-intersecting triples $W_1,W_2,W_3$ of finite binary walks having starting points $(-1,1)$, $(0,0)$ and $(0,-1)$,
same numbers of North-steps, and with $\mathrm{length}(W_1)=\mathrm{length}(W_2)=\mathrm{length}(W_3)-1$.  
As described in~\cite{BoBoFu09} (the original bijection, between Baxter permutations
and triples of paths, is due to Dulucq and Guibert~\cite{dulucq-guibert-baxter}) there is a bijection
between the class $\cO$ of plane bipolar orientations and the class $\cT$, 
with following parameter-correspondence (see also Figure~\ref{fig:bije}(a)):

\begin{center}
\begin{tabular}{rcl}
 $n$ edges & $\leftrightarrow$ & $\mathrm{length}(W_1)=n-1$, \\
$k$ non-pole vertices & $\leftrightarrow$ & $W_1$ has $k$ East-steps,\\ 
$i+1$ sink-edges & $\leftrightarrow$ & $W_3$ ends with an East-step followed by $i$ North-steps. 
 \end{tabular}
 \end{center}

Now define $\cE$ as the class of 4-tuples $(W_1,W_2,W_3,S)$, where $(W_1,W_2,W_3)$
is a non-intersecting triple of binary walks starting respectively from the points $(-1,1)$, $(0,0)$ and $(0,-1)$; and $S$ is a sequence of non-negative integers such that:
\begin{itemize}
\item
The walks $W_1$ and $W_2$ have same lengths and same numbers of East-steps.
\item
Denoting by $(x_2,y_2)$ and $(x_3,y_3)$ the coordinates of the respective end-points of
$W_2$ and $W_3$, we have $x_3\geq x_2$ and $x_2+y_2\geq x_3+y_3$.
Let $a:=x_3-x_2$ and $b:=x_2+y_2-x_3-y_3$.
\item
The sequence $S$ is made of $a+1$ non-negative integers that add up to $b$.
\end{itemize}

Let $O$ be a decorated plane bipolar orientation with sink-degree $i+1$.      
Let $(W_1,W_2,W_3)\in\cT$ be the 
  triple of walks associated with $O$ (without the decorations), 
and $(W,S)$ the sink-code for $O$. 
Define $W_3'$ as $W_3$ where the suffix $\mathrm{East}\ \mathrm{North}^i$ is replaced by $W$, see Figure~\ref{fig:bije}. 
Then it is easily checked that $E:=(W_1,W_2,W_3',S)\in\cE$.
Note that one can recover $W_3$ from $E$ ($W_3$ is the unique binary walk ending at $(x_2+1,y_2-1)$ and equal to $W_3'$ in the area $\{x\leq x_2\}$), hence
one can recover $O$. Moreover if $O$ has $p$ marked sink-edges and $q$ marked sink-corners,
then with the notations above, $b=p$ and $a=q+1-p$. 

Overall, we obtain a bijection between the class $\cD$ of decorated
plane bipolar orientations and the class $\cE$, with the following parameter-correspondence:

\begin{center}
\begin{tabular}{rcl}
 $m$ edges  & $\leftrightarrow$ &  $W_1$ and $W_2$  have length $m-1$,\\ 
 $k$ non-pole vertices & $\leftrightarrow$ &  $W_1$ and $W_2$ have $k$ East-steps,\\
 $p$ marked sink-edges & $\leftrightarrow$ &  $W_3$ has length $m-p$,\\
 $q$ marked sink-corners &  $\leftrightarrow$ & $W_3$ has $k+q-p+1$ East-steps. 
 \end{tabular}
 \end{center}

Composing the bijection between $\cI$ and $\cD$ with the bijection between $\cD$ and $\cE$
we finally obtain (taking $m=n+p$ and $q=p+r-1$ in the correspondence above):

\begin{Theorem}\label{theo:bij}
There is a bijection between involutive Baxter permutations with  
$2n$ non-fixed points, $2k$ descents not crossing the diagonal, $p$ fixed points, $r$ descents
crossing the diagonal; and $4$-tuples of the form $(W_1,W_2,W_3,S)$
where $(W_1,W_2,W_3)$ is a non-intersecting triple of binary walks
with starting points $(-1,1)$, $(0,0)$, $(0,-1)$, end-points $(k-1, n+p-k)$, $(k,n+p-k-1)$, 
$(k+r,n-k-r-1)$, and where $S$ is a sequence of $r+1$ nonnegative numbers adding up to $p$.
\end{Theorem}

\section{Counting}
Let $(A_1,A_2,A_3)$ and $(B_1,B_2,B_3)$ be the starting points
and end-points in Theorem~\ref{theo:bij}. By the Lindstr\"om-Gessel-Viennot
lemma~\cite{gessel-viennot}, the number of non-intersecting triples
of binary walks with starting points $A_1,A_2,A_3$
and end-points $B_1,B_2,B_3$ is the determinant of the $3\times 3$ matrix
$(m_{i,j})$, with $m_{i,j}$ the number of binary walks from $A_i$ to $B_j$
(note that each entry $m_{i,j}$ is an explicit binomial coefficient).
The number of choices for the sequence $S$ in Theorem~\ref{theo:bij}
is clearly equal to $\binom{p+r}{r}$. We obtain (taking out of the 
determinant a common binomial factor for each column):

\begin{Theorem}[multivariate enumeration formula]\label{theo:multi}
For $n>0$, and $k$, $p$, $r$ nonnegative integers, the number $a_{n,k,p,r}$ of involutive Baxter permutations with  
$2n$ non-fixed points, $2k$ descents not crossing the diagonal, $p$ fixed points, 
$r$ descents
crossing the diagonal is given by
\beq
a_{n,k,p,r}=\frac{\displaystyle\binom{p+r}{r}\binom{n+p-1}{k}^2\binom{n}{t}}{nq^2(q\!+\!1)(k\!+\!1)(t\!+\!1)}\cdot 
\left|\begin{array}{ccc}
q(q\!+\!1)&q(q\!-\!1)&s(s\!-\!1)\\[.2cm]
k(q\!+\!1)&(k\!+\!1)q&s(t\!+\!1)\\[.2cm]
k(k\!-\!1)&k(k\!+\!1)&t(t\!+\!1) 
\end{array}\right|
\eeq
where $q:=n+p-k$, $s:=n-k-r$, $t:=k+r$.
\end{Theorem}
An equivalent multivariate formula for $a_{n,k,0,r}$ has been
obtained by Bousquet-M\'elou~\cite{bo10} using the ``obstinate'' kernel method.  
Note that, by the correspondence of Section~\ref{sec:baxter_decorated_plane}, 
the number $a_{n,k,0,r}$ counts monosource orientations (without taking sailing sinks into account)
 with $n$ edges, $r$ sinks, and $k$ non-extremal vertices. 

We now prove that the number $b_n$ of involutive
Baxter permutations with no fixed point and $2n$ elements
satisfies $b_n=\frac{3\cdot 2^{n-1}}{(n+1)(n+2)}\binom{2n}{n}$   
(even though $b_n=\sum_{k,r}a_{n,k,0,r}$,  our proof does not exploit 
 the multivariate formula of Theorem~\ref{theo:multi}). 
Let $\cF_n$ be the set of involutive Baxter permutations of size $2n$ with no fixed point. 
In this part it is convenient to  rotate the encoding triples of paths by $45$ degrees counter-clockwise, 
to delete the first (East-) step in the third path of the triple,  
and to rescale by $\sqrt{2}$. This way, the binary walks considered in the previous section become 
paths having steps $(-1/2,+1/2)$ 
(annotated $\nwarrow$) or $(+1/2,+1/2)$ (annotated $\nearrow$)~\footnote{We talk about \emph{walks} when steps are \{East,North\} and \emph{paths} when steps are $\{\nwarrow,\nearrow\}$.}, 
and the bijection of Theorem~\ref{theo:bij} specializes as follows:

\begin{claim}
For $n\geq 1$, $\cF_n$ is in bijection with the set $\cR_n$ 
of non-intersecting triples of paths, each with $n-1$ steps either $\nwarrow$ or $\nearrow$,
with starting points $(-1,0)$, $(0,0)$, $(1,0)$ and such that the endpoints
of the first two paths are at distance
$1$ on the line $\{y=\tfrac{n-1}{2}\}$. 
\end{claim}
\begin{proof}
By Theorem~\ref{theo:bij} $\cF_n$ is in bijection with the set
of non-intersecting triples of finite binary walks $(W_1,W_2,W_3)$ with starting points $(-1,1)$, $(0,0)$, and $(0,-1)$, and end-points on the line
$\{x+y=n-1\}$ such that the end-points of $W_1$ and $W_2$ are at distance $\sqrt{2}$. Since the $3$ walks do not intersect, 
the first step of the third walk is always an East-step, hence can be removed
with no loss of information. This way the $3$ walks start
on the line $\{x+y=0\}$ with $W_1$ at distance $\sqrt{2}$ from $W_2$ itself
at distance $\sqrt{2}$ from $W_3$. Hence, rotating the figure 
counter-clockwise by $\pi/4$
and rescaling by $\sqrt{2}$, one has a triple of paths
in $\cR_n$. 
\end{proof}

To enumerate the triples in $\cR_n$, we inject $\cR_n$ in the bigger
set $\cU_n$ defined the same way, except that 
the third path is allowed to intersect the two other paths. Let $\cS_n$ be the
subset of objects in $\cU_n$ where the third path meets
the second path (and possibly also the first path); so we have
$\cU_n=\cR_n+\cS_n$. Let $u_n$, $r_n$, $s_n$ be the cardinalities of
  $\cU_n$, $\cR_n$, and $\cS_n$, respectively.  
We have $u_n=r_n+s_n$, so that 
$$
b_n=r_n=u_n-s_n.
$$

\begin{figure}[htb!]
\begin{center}
\includegraphics[width=12cm]{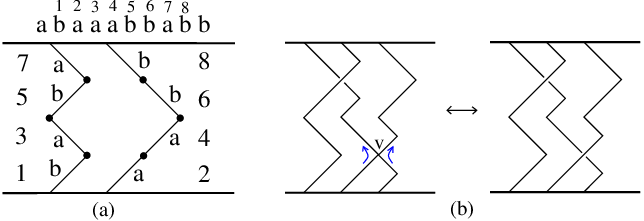}
\caption{(a) A non-intersecting pair of paths with starting points
and end-points at distance $1$ can be encoded by a Dyck word, 
(b) bijection
between $\cV_n$ and $\cS_n+\mathrm{mir}(\cS_n)$.}
\label{fig:paths}
\end{center}
\end{figure}

We now state a basic lemma~\cite{Lev:59} 
which we will use to obtain formulas for $u_n$ and $s_n$ 
(the case $k=1$, which gives the Catalan numbers,  
is illustrated in Figure~\ref{fig:paths}(a)): 

\begin{Lemma}[folklore]\label{lem:folklore}
For $n$ and $k$ positive integers, let $a_n^{(k)}$ be the number of non-intersecting 
pairs $(P_1,P_2)$ of paths (counted up to horizontal translation) each with 
$n-1$ steps either $\nwarrow$ or $\nearrow$, starting at distance $k$
on the line $\{y=0\}$ and ending at distance $1$ on the line $\{y=(n-1)/2\}$.
Then
$$
a_n^{(k)}=\frac{2k(2n-1)!}{(n-k)!(n+k)!}.
$$
\end{Lemma}

The lemma directly yields a formula for $u_n$; indeed the first two
paths of a triple in $\cU_n$ are non-intersecting,  
start at distance $1$ on $\{y=0\}$ and end at distance $1$ on $\{y=\tfrac{n-1}{2}\}$, and the third path is unconstrained. Hence
$$
u_n=2^{n-1}a_n^{(1)}=2^{n-1}\frac{(2n)!}{n!(n+1)!}.
$$
To obtain a formula for $s_n$,
we ``double'' the set $\cS_n$.
For a triple $\gamma\in\cS_n$, the \emph{mirror} of $\gamma$ is obtained by applying a vertical mirror (reflexion
according to a vertical line, up to a global translation) to $\gamma$. 
Denote by $\mathrm{mir}(\cS_n)$ the set of mirrors of triples in $\cS$.
Note that  $\cS_n\cap\mathrm{mir}(\cS_n)$ is empty (indeed, in $\cS_n$, the middle-starting
path intersects only 
the right-starting path; whereas in $\mathrm{mir}(\cS_n)$ the middle-starting
path intersects only  the 
left-starting path). 
We now establish a bijection between $\cS_n+\mathrm{mir}(\cS_n)$ and a set easy
to count. Let $\cV_n$ be the set of triples of paths, each with $n-1$ steps
either $\nwarrow$ or $\nearrow$, with starting points $(-1,0)$, $(0,0)$, $(1,0)$,
and such that the left-starting path and the right-starting path do not intersect and end
at distance $1$ on the line $\{y=\tfrac{n-1}{2}\}$.

\begin{claim}\label{claim:Vn}
The set $\cV_n$ is in bijection with $\cS_n+\mathrm{mir}(\cS_n)$. 
\end{claim}
\begin{proof}
The bijection relies on a simple argument akin to the Gessel-Viennot lemma. First we describe
the mapping from $\cV_n$ to $\cS_n+\mathrm{mir}(\cS_n)$. 
Consider a triple of paths from $\cV_n$, and 
denote by $P_{\ell}$, $P_m$, $P_r$ the paths starting from $(-1,0)$, $(0,0)$, and $(1,0)$,
respectively. 
Since the end-points of $P_{\ell}$ and $P_r$ are at distance $1$ (i.e., consecutive) on $\{y=\tfrac{n-1}{2}\}$,
the path $P_m$ has to intersect $P_{\ell}\cup P_r$. 
Let $v$ be the first intersection of $P_m$ with $P_{\ell}\cup P_r$ (note that $v$ can not be on both $P_{\ell}$ and on $P_r$, 
see Figure~\ref{fig:paths}(b)).  
If $v\in P_r$ we exchange the parts of $P_r$ and $P_m$ after 
$v$, this yields a triple of paths in $\cS$ (see Figure~\ref{fig:paths}(b));  
if $v\in P_{\ell}$ we exchange the parts of $P_{\ell}$ and $P_m$ after 
$v$, this yields a triple of paths in $\mathrm{mir}(\cS)$.
It is now straightforward to get 
the inverse mapping, from $\cS_n+\mathrm{mir}(\cS_n)$ to $\cV_n$. 
Let $\gamma\in\cS_n+\mathrm{mir}(\cS_n)$, and 
denote again by $P_{\ell}$, $P_m$, and $P_r$ the paths starting from $(-1,0)$, $(0,0)$, and $(1,0)$,
respectively. 
If $\gamma\in\cS$ let $v$ be the first intersection of $P_m$ and $P_r$; exchange 
the portions of $P_m$ and $P_r$ after $v$. 
If $\gamma\in\mathrm{mir}(\cS)$ let $v$ be the first intersection of $P_{\ell}$ and $P_m$; exchange the portions of $P_{\ell}$ and $P_m$ after $v$. 
\end{proof}
Let $v_n$ be the cardinality of $\cV_n$. 
Claim~\ref{claim:Vn} implies that $v_n=2s_n$. Now $v_n$ is easy to obtain. Indeed, for a triple in $\cV_n$,    
the left-starting path and right-starting path start at distance $2$ on $\{y=0\}$, 
end at distance $1$ on $\{y=\frac{n-1}{2}\}$, 
and are non-intersecting; and the middle-starting path is unconstrained.
Hence, with the notation of Lemma~\ref{lem:folklore}, 
$$
v_n=2^{n-1}a_{n}^{(2)}=2^{n+1}\frac{(2n-1)!}{(n-2)!(n+2)!},
$$
so that $s_n=v_n/2= 2^{n}\frac{(2n-1)!}{(n-2)!(n+2)!}$.
From $b_n=u_n-s_n$ and the expressions of $u_n$ and $s_n$ we obtain: 

\begin{Theorem}[univariate formula, recovers~\cite{bo10} in a bijective way]\label{theo:uni}
The number $b_n$ of involutive Baxter permutations with no fixed
point and $2n$ elements is
\beq
b_n=\frac{3\cdot 2^{n-1}}{(n+1)(n+2)}\binom{2n}{n}.
\eeq
\end{Theorem}
Note that, by the correspondence of Section~\ref{sec:baxter_decorated_plane}, $b_n$ is the number of monosource orientations (without
taking sailing sinks into account) with $n$ edges. 

\vspace{.4cm}

\noindent\emph{Acknowledgement.} I am very grateful to Mireille Bousquet-M\'elou
and Nicolas Bonichon for helfpul discussions and corrections on a preliminary version.
I also thank Mireille 
for communicating to me the counting formulas for involutive Baxter permutations. 

\bibliographystyle{plain}
\bibliography{baxter.bib}

\end{document}